\documentclass[reqno,twoside]{amsart}
      \usepackage{amsmath}  \usepackage{amsthm}
      \usepackage{amssymb}  \usepackage{latexsym}
  \newtheorem{theorem}{Theorem}[section]
  \newtheorem{lemma}[theorem]{Lemma}
  \newtheorem{proposition}[theorem]{Proposition}
  \newtheorem{corollary}[theorem]{Corollary}
  \theoremstyle{definition}
  \newtheorem{definition}[theorem]{Definition}
  \newtheorem{example}[theorem]{Example}

\begin{document}

\title{Fuzzy $n$-ary groups as a generalization of Rosenfeld's fuzzy groups}
\author{B. Davvaz and  Wies{\l}aw A. Dudek}

\address{B. Davvaz,
Department of Mathematics, Yazd University, Yazd, Iran}
\email{davvaz@yazduni.ac.ir}

\address{W.A. Dudek,
Institute of Mathematics and Computer Science,  Wroc{\l}aw
University of Technology, Wybrze\.ze Wyspia\'nskiego 27, 50-370
Wroc{\l}aw, Poland}
\email{dudek@im.pwr.wroc.pl}

 {

\begin{abstract}
The notion of an $n$-ary group is a natural generalization of the
notion of a group and has many applications in different branches.
In this paper, the notion of (normal) fuzzy $n$-ary subgroup of an
$n$-ary group is introduced and some related properties are
investigated. Characterizations of fuzzy $n$-ary subgroups are
given.
\end{abstract}
\maketitle

\section{Preliminaries}
A nonempty set $G$ together with one $n$-ary operation $f:G^n
\longrightarrow G$, where $n\geq 2$, is called an {$n$-ary
groupoid} and is denoted by $(G,f)$. According to the general
convention used in the theory of $n$-ary groupoids the sequence of
elements $x_i,x_{i+1},\ldots,x_j$ is denoted by $x_i^j$. In the
case $j<i$ it is the empty symbol. If
$x_{i+1}=x_{i+2}=\ldots=x_{i+t}=x$, then instead of
$x_{i+1}^{i+t}$ we write $\stackrel{(t)}{x}$. In this convention
$f(x_1,\ldots,x_n)= f(x_1^n)$ and
 \[
 f(x_1,\ldots,x_i,\underbrace{x,\ldots,x}_{t},x_{i+t+1},\ldots,x_n)=
 f(x_1^i,\stackrel{(t)}{x},x_{i+t+1}^n) .
 \]

An $n$-ary groupoid $(G,f)$ is called {\it $(i,j)$-associative} if
 \begin{equation}
f(x_1^{i-1},f(x_i^{n+i-1}),x_{n+i}^{2n-1})=
f(x_1^{j-1},f(x_j^{n+j-1}),x_{n+j}^{2n-1})\label{assolaw}
 \end{equation}
holds for all $x_1,\ldots,x_{2n-1}\in G$. If this identity holds
for all $1\leqslant i<j\leqslant n$, then we say that the
operation $f$ is {\it associative} and $(G,f)$ is called an {\it
$n$-ary semigroup}. It is clear that an $n$-ary groupoid is
associative if and only if it is $(1,j)$-associative for all
$j=2,\ldots,n$. In the binary case (i.e. for $n=2$) it is a usual
semigroup.

If for all $x_0,x_1,\ldots,x_n\in G$ and fixed
$i\in\{1,\ldots,n\}$ there exists an element $z\in G$ such that
\begin{equation}                                                          \label{solv}
f(x_1^{i-1},z,x_{i+1}^n)=x_0 ,
\end{equation}
then we say that this equation is {\it $i$-solvable} or {\it
solvable at the place $i$}. If this solution is unique, then we
say that (\ref{solv}) is {\it uniquely $i$-solvable}.

An $n$-ary groupoid $(G,f)$ uniquely solvable for all
$i=1,\ldots,n$ is called an {\it $n$-ary quasigroup}. An
associative $n$-ary quasigroup is called an {\it $n$-ary group}.
It is clear that for $n=2$ we obtain a usual group.

Note by the way that in many papers $n$-ary semigroups ($n$-ary
groups) are called $n$-semigroups ($n$-groups, respectively).
Moreover, in many papers, where the arity of the basic operation
does not play a crucial role, we can find the term a {\it polyadic
semigroup} ({\it polyadic group}) (cf. \cite{Post}).

Now such and similar $n$-ary systems have many applications in
different branches. For example, in the theory of automata
\cite{Busse} $n$-ary semigroups and $n$-ary groups are used, some
$n$-ary groupoids are applied in the theory of quantum groups
\cite{Nik}. Different applications of ternary structures in
physics are described by R. Kerner in \cite{Ker}. In physics there
are used also such structures as $n$-ary Filippov algebras (see
\cite{Poj}) and $n$-Lie algebras (see \cite{Vai}).

The idea of investigations of such groups seems to be going back to
E. Kasner's lecture \cite{kasner} at the fifty-third annual meeting
of the American Association for the Advancement of Science in 1904.
But the first important paper concerning the theory of $n$-ary
groups was written (under inspiration of Emmy Noether) by W.
D\"ornte in 1928 (see \cite{Dor}). In this paper D\"ornte observed
that any $n$-ary groupoid $(G,f)$ of the form $\,f(x_1^n)=x_1\circ
x_2\circ\ldots\circ x_n\circ b$, where $(G,\circ)$ is a group and
$b$ belongs to the center of this group, is an $n$-ary group but for
every $n>2$ there are $n$-ary groups which are not of this form. In
the first case we say that an $n$-ary group $(G,f)$ is {\it
$b$-derived} (or {\it derived} if $b$ is the identity of
$(G,\circ)$) from the group $(G,\circ )$, in the second -- {\it
irreducible}. Moreover, in some $n$-ary groups there exists an
element $e$ (called an {\it $n$-ary neutral element}) such that
 \begin{equation}                                            \label{n-id}
 f(\stackrel{(i-1)}{e},x,\stackrel{(n-i)}{e})=x
 \end{equation}
holds for all $x\in G$ and for all $i=1,\ldots,n$. It is
interesting that each $n$-ary group $(G,f)$ containing a neutral
element are derived from a binary group $(G,\circ)$, where $x\circ
y=f(x,\stackrel{(n-2)}{e},y)$ (cf. \cite{Dor}). On the other hand,
there are $n$-ary groups with two, three and more neutral
elements. All $n$-ary groups with this property are derived from
the commutative group of the exponent $k|(n-1)$.

It is worthwhile to note that in the definition of an $n$-ary
group, under the assumption of the associativity of $f$, it
suffices only to postulate the existence of a solution of
(\ref{solv}) at the places $i=1$ and $i=n$ or at one place $i$
other than $1$ and $n$. Then one can prove the uniqueness of the
solution of (\ref{solv}) for all $i=1,\ldots,n$ (cf. \cite{Post},
p. $213^{17}$). Some other definitions of $n$-ary groups one can
find in \cite{Rem} and \cite{DGG}.

In an $n$-ary group the role of the inverse element plays the
so-called {\it skew element}, i.e., an element $\overline{x}$ such
that
 \begin{equation}\label{skew}
f(\stackrel{(n-1)}{x},\overline{x})=x.
 \end{equation}
It is uniquely determined, but $\overline{x}=\overline{y}$ do not
implies $x=y$, in general. Moreover, there are $n$-ary groups in
which one element is skew to all (cf. \cite{D90}). So, in general,
the skew element to $\overline{x}$ is not equal to $x$, but in
ternary $(n=3)$ groups we have $\overline{\overline{x}}=x$ for all
$x\in G$. For some elements of $n$-ary groups we have
$\overline{x}=x$. Such elements are called {\it idempotents}. An
$n$-ary group in which elements are idempotents is called {\it
idempotent}. There are $n$-ary groups without idempotents. A
simple example of $n$-ary groups with only one idempotent are {\it
unipotent} $n$-ary groups described in \cite{Uni}. In these groups
there exists an element $\theta$ such that $f(x_1^n)=\theta$ holds
for all $x_1^n\in G$.

Note that in all $n$-ary groups the following two identities
\begin{equation}\label{e-dor}
f(y,\stackrel{(i-2)}{x},\overline{x},\stackrel{(n-i)}{x})=
f(\stackrel{(n-j)}{x},\overline{x},\stackrel{(j-2)}{x},y)=y
\end{equation}
are satisfied for all $2\leq i,j\leq n$ (cf. \cite{Dor}).

A nonempty subset $H$ of an $n$-ary group $(G,f)$ is an {\it
$n$-ary subgroup} if $(H,f)$ is an $n$-ary group, i.e., if it is
closed under the operation $f$ and $x\in H$ implies
$\overline{x}\in H$ (cf. \cite{Dor}). The intersection of two
subgroups may be the empty set. Moreover, there are $n$-ary groups
which are set theoretic union of disjoint isomorphic subgroups
(cf. for example \cite{auto}).

Fixing in an $n$-ary operation $f$, where $n\geq 3$, the elements
$a_2^{n-1}$ we obtain the new binary operation $x\diamond
y=f(x,a_2^{n-2},y)$. If $(G,f)$ is an $n$-ary group then
$(G,\diamond)$ is a group. Choosing different elements $a_2^{n-2}$
we obtain different groups. All these groups are isomorphic
\cite{DM}. So, we can consider only groups of the form
$ret_a(G,f)=(G,\circ)$, where $x\circ y=
f(x,\stackrel{(n-2)}{a},y)$. In this group $e=\overline{a}$,
$x^{-1}=f(\overline{a},\stackrel{(n-3)}{x},\overline{x},\overline{a})$.

Subgroups of $(G,f)$ are not subgroups of $(G,\circ)$, in general.

In the theory of $n$-ary groups a very important role plays the
following theorem firstly proved by M. Hossz\'u \cite{Hos} (see
also \cite{DM1}).

\begin{theorem}\label{T1.1}
For any $n$-ary group $(G,f)$ there exist a group $(G,\circ)$, its
automorphism $\varphi$ and an element $b\in G$ such that
 \begin{equation}\label{ho}
f(x_1^n)=x_1\circ\varphi(x_2)\circ\varphi^2(x_3)\circ\ldots\circ\varphi^{n-1}(x_n)\circ
b
 \end{equation}
 holds for all $x_1^n\in G$.
\end{theorem}

One can proved (see for example \cite{DM1}) that in this theorem
$(G,\circ)=ret_a(G,f)$,
$\varphi(x)=f(\overline{a},x,\stackrel{(n-2)}{a})$,
$b=f(\overline{a},\ldots,\overline{a})$, where $a$ is an arbitrary
element of $G$. The above representation is unique up to
isomorphism.

Since, as it is not difficult to see,
\[
\varphi^{n-1}(x)\circ b=     b\circ x
\]
the identity (\ref{ho}) can be written in more useful form
 \begin{equation}\label{ho-2}
f(x_1^n)=x_1\circ\varphi(x_2)\circ\varphi^2(x_3)\circ\ldots\circ\varphi^{n-2}(x_n)\circ
b\circ x_n.
 \end{equation}

\section{Fuzzy $n$-ary subgroups}

Any function $\mu :G\longrightarrow [0,1]$ is called a {\it fuzzy
subset} of $G$. The set of all values of $\mu$ is denoted by
Im$(\mu)$. If for every $S\subseteq G$, there exists $x_0\in S$
such that $\mu(x_0)=\sup\{\mu(x)\,|\,x\in S\}$ then we say that
$\mu$ has {\it sup-property}.

For usual groups A. Rosenfeld defined \cite{Ro} fuzzy subgroups in
the following way:
\begin{definition}
A fuzzy subset $\mu$ defined on a group $(G,\cdot)$ is called a
{\it fuzzy subgroup } if
 \begin{enumerate}
\item[1)] $\mu(xy)\geq\min\{\mu(x),\mu(y)\}$,
\item[2)] $\mu(x^{-1})\geq\mu(x)$
\end{enumerate}
holds for all $x,y\in G$.
\end{definition}
In fact we have $\mu(x^{-1})=\mu(x)$ because $(x^{-1})^{-1}=x$ for
every $x\in G$. Moreover, from the above definition we can deduce
that $\mu(e)\geq\mu(x)$ for every $x\in G$.

\begin{proposition} {\rm \cite{Ro}}\label{P2.2}
A fuzzy subset $\mu$ on a group $(G,\cdot)$ is a fuzzy subgroup if
and only if each nonempty level subset $\mu_t=\{x\in
G\,|\,\mu(x)\geq t\}$ is a subgroup of $(G,\cdot)$.
\end{proposition}

The above definition can be extended to $n$-ary case in the
following way (cf. \cite{QRS8}):

\begin{definition}\label{D2.3} Let $(G,f)$ be an $n$-ary group. A
fuzzy subset of $G$ is called a {\it fuzzy $n$-ary subgroup} of
$(G,f)$ if the following axioms hold:
\begin{itemize}
\item[$(i)$] \ $\mu(f(x_1^n))\geq\min\{\mu(x_1),\ldots,\mu(x_n)\}$ for all $x_1^n\in G$,
\item[$(ii)$] \ $\mu(\overline{x})\geq\mu(x)$ for all $x\in G$.
\end{itemize}
\end{definition}

Note that for $n=3$ the second condition of the Definition
\ref{D2.3} can be replaced by the condition
\begin{enumerate}
\item[$(iii)$] $\mu(\overline{x})=\mu(x)$ for all $x\in G$,
\end{enumerate}
because in this case $n=3$ we have $\overline{\overline{x}}=x$ for
all $x\in G$ (cf. \cite{Dor}). These two conditions are equivalent
for all $n$-ary groups in which for every $x\in G$ there exists a
natural number $k$ such that $\bar{x}^{(k)}=x$, where
$\bar{x}^{(k)}$ denotes the element skew to $\bar{x}^{(k-1)}$ and
$\bar{x}^{(0)}=x$. But, as it was observed in \cite{QRS8}, there
are fuzzy $n$-ary subgroups in which $\mu(\bar{x})>\mu(x)$ for
some $x\in G$.

\begin{example}\label{E2.4} Let $(\mathbb{Z}_4,f)$ be a $4$-ary group
derived from the additive group $\mathbb{Z}_4$. It is not difficult
to see that the map $\mu$ defined by $\mu(0)=1$ and $\mu(x)=0.5$ for
all $x\ne 0$ is a fuzzy $4$-ary subgroup in which for $x=2$ we have
$\overline{x}=0$ and $\mu(\bar{x})>\mu(x)$.
\end{example}

\begin{proposition}\label{P2.5} Any $n$-ary subgroup of $(G,f)$ can be realized as
a level subset of some fuzzy $n$-ary subgroup of $G$.
\end{proposition}
\begin{proof} Let $H$ be an $n$-ary subgroup of a given
$n$-ary group $(G,f)$ and let $\mu_H$ be a fuzzy subset of $G$
defined by
\[
\mu_H(x)= \left\{\begin{array}{ll} t & \text {if } \ x\in H\\
s &  \text {if} \ \ x \not\in H
\end{array}
\right.
\]
where $0\leq s<t\leq 1$ is fixed. It is not difficult to see that
$\mu$ is a fuzzy $n$-ary subgroup of $G$ such that $\mu_t =H$.
\end{proof}

\begin{corollary}\label{C2.6}
The characteristic function of a nonempty subset of an $n$-ary
group $(G,f)$ is a fuzzy $n$-ary subgroup of $G$ if and only if
$A$ is an $n$-ary subgroup of $G$.
\end{corollary}

\begin{theorem}\label{T2.7}
A fuzzy subset $\mu$ on an $n$-ary group $(G,f)$ is a fuzzy
$n$-ary subgroup if and only if each its nonempty level subset is
an $n$-ary subgroup of $(G,\cdot)$.
\end{theorem}
\begin{proof} Let $\mu$ be a fuzzy $n$-ary subgroup of an $n$-ary group $(G,f)$. If
$x_1^n\in\mu_t$ for some $t\in [0,1]$, then $\mu(x_i)\geq t$ for
all $i=1,2,\ldots,n$. Thus
\[
\mu(f(x_1^n))\geq\min\{\mu(x_1),\ldots,\mu(x_n)\}\geq t,
\]
which implies $f(x_1^n)\in\mu_t$. Moreover, for $x\in\mu_t$ from
$\mu(\overline{x})\geq\mu(x)\geq t$ it follows
$\overline{x}\in\mu_t$. So, $\mu_t$ is an $n$-ary subgroup of
$(G,f)$.

Conversely, assume that every nonempty level subset $\mu_t$ is an
$n$-ary subgroup of $(G,f)$. Let
$t_0=\min\{\mu(x_1),\ldots,\mu(x_n)\}$ for some $x_1^n\in G$. Then
obviously $x_1^n\in\mu_{t_0}$, consequently,
$f(x_1^n)\in\mu_{t_0}$. Thus
\[
\mu(f(x_1^n))\geq t_0=\min\{\mu(x_1),\ldots,\mu(x_n)\}.
\]
Now let $x\in\mu_t$. Then $\mu(x)=t_0\geq t$, i.e.,
$x\in\mu_{t_0}$. Since, by the assumption, every nonempty level
set of $\mu$ is an $n$-ary subgroup, $\overline{x}\in\mu_{t_0}$.
Whence $\mu(\overline{x})\geq t_0=\mu(x)$.

In this way the conditions of Definition \ref{D2.3} are verified.
This completes the proof.
\end{proof}

Using this theorem we can prove the another characterization of
fuzzy $n$-ary subgroups.

\begin{theorem}\label{T2.8}
A fuzzy subset $\mu$ on an $n$-ary group $(G,f)$ is a fuzzy
$n$-ary subgroup if and only if for all $i=1,2,\ldots,n$ and all
$x_1^n\in G$ it satisfies the following two conditions
\begin{itemize}
\item[$(i)$] \ $\mu(f(x_1^n))\geq\min\{\mu(x_1),\ldots,\mu(x_n)\}$,
\item[$(ii)$] \ $\mu(x_i)\geq\min\{\mu(x_1),\ldots,\mu(x_{i-1}),
\mu(f(x_1^n)),\mu(x_{i-1}),\ldots,\mu(x_n)\}$.
\end{itemize}
\end{theorem}
\begin{proof}
Assume that $\mu$ is a fuzzy $n$-ary subgroup of $(G,f)$.
Similarly as in the proof of Theorem~~\ref{T2.7} we can prove that
each nonempty level subset $\mu_t$ is closed under the operation
$f$, i.e., $x_1^n\in\mu_t$ implies $f(x_1^n)\in\mu_t$.

Now let $x_0,x_1^{i-1},x_{i+1}^n$, where
$x_0=f(x_1^{i-1},z,x_{i+1}^n)$ for some $i=1,2,\ldots,n$ and $z\in
G$, be in $\mu_t$. Then, according to $(ii)$, $\mu(z)\geq t$,
which proves $z\in\mu_t$. So, the equation (\ref{solv}) has a
solution $z\in\mu_t$. This means that each nonempty $\mu_t$ is an
$n$-ary subgroup.

Conversely, if all nonempty level subsets of $\mu$ are $n$-ary
subgroups, than, similarly as in the previous proof, we can see
that the condition $(i)$ is satisfied. Moreover, if for $x_1^n\in
G$ we have
\[
t=\min\{\mu(x_1),\ldots,\mu(x_{i-1}),
\mu(f(x_1^n)),\mu(x_{i-1}),\ldots,\mu(x_n)\},
\]
then $x_1^{i-1},x_{i+1}^n,f(x_1^n)\in\mu_t$. Whence, according to
the definition of an $n$-ary group, we conclude $x_i\in\mu_t$.
Thus $\mu(x_i)\geq t$. This proves $(ii)$.
\end{proof}

\begin{corollary}\label{C2.9}
A fuzzy subset $\mu$ defined on a group $(G,\cdot)$ is a {\it
fuzzy subgroup } if and only if
 \begin{enumerate}
\item[1)] $\mu(xy)\geq\min\{\mu(x),\mu(y)\}$,
\item[2)] $\mu(x)\geq\min\{\mu(y),\mu(xy)\}$,
\item[3)] $\mu(y)\geq\min\{\mu(x),\mu(xy)\}$
\end{enumerate}
holds for all $x,y\in G$.
\end{corollary}

\begin{theorem}\label{T2.10}
Let $\mu$ be a fuzzy $n$-ary subgroup of $(G,f)$. If there exists
an element $a\in G$ such that $\mu(a)\geq\mu(x)$ for every $x\in
G$, then $\mu$ is a fuzzy subgroup of a group $ret_a(G,f)$.
\end{theorem}
\begin{proof}
Indeed, \[\mu(x\circ
y)=\mu(f(x,\stackrel{(n-2)}{a},y)\geq\min\{\mu(x),\mu(a),\mu(y)\}
=\min\{\mu(x),\mu(y)\}
\]
and
\[
\mu(x^{-1})=\mu(f(\overline{a},\stackrel{(n-3)}{x},\overline{x},\overline{a}))\geq
\min\{\mu(x),\mu(\overline{x}),\mu(a),\mu(\overline{a})\}=\mu(x),
\]
which completes the proof.
\end{proof}

The assumption that $\mu(a)\geq\mu(x)$ cannot be omitted.

\begin{example}\label{E2.11}
Let $(\mathbb{Z}_4,f)$ be a ternary group from Example
~\ref{E2.4}. Then a fuzzy set $\mu$ defined by $\mu(0)=1$,
$\mu(2)=0.5$, $\mu(1)=\mu(3)=0.3$ is a fuzzy ternary subgroup of
$(\mathbb{Z}_4,f)$. (This fact follows also from our
Theorem~\ref{T3.5} because $\{0\}$ and $\{0,2\}$ are subgroups of
$(\mathbb{Z}_4,f)$.) For $ret_1(Z_4,f)$ we have $\mu(2\circ
2)=\mu(f(2,1,2))=\mu(1)=0.3<\min\{\mu(2),\mu(2)\}=0.5$. So, the
assumption $\mu(a)\geq\mu(x)$ cannot be omitted.
\end{example}

\begin{theorem}\label{T2.12}
Let $(G,f)$ ba an $n$-ary group. If $\mu$ is a fuzzy subgroup of a
group $ret_a(G,f)$ and $\mu(a)\geq\mu(x)$ for all $x\in G$, then
$\mu$ is a fuzzy $n$-ary group of $(G,f)$.
\end{theorem}
\begin{proof}
According to Theorem~\ref{T1.1} any $n$-ary group can be presented
in the form (\ref{ho}), where $(G,\circ )=ret_a(G,f)$,
$\varphi(x)=f(\overline{a},x,\stackrel{(n-2)}{a})$ and
$b=f(\overline{a},\ldots,\overline{a})$. Obviously,
\[
\mu(\varphi(x))=\mu(f(\overline{a},x,\stackrel{(n-2)}{a}))\geq
\min\{\mu(\overline{a}),\mu(x),\mu(a)\}=\mu(x),
\]
\[
\mu(\varphi^2(x))=\mu(f(\overline{a},\varphi(x),\stackrel{(n-2)}{a}))\geq
\min\{\mu(\overline{a}),\mu(\varphi(x)),\mu(a)\}=\mu(\varphi(x))\geq\mu(x).
\]
Consequently, $\mu(\varphi^k(x))\geq\mu(x)$ for all $x\in G$ and
$k\in\mathbb{N}$.

Similarly
\[
\mu(b)=\mu(f(\overline{a},\ldots,\overline{a}))\geq\mu(\overline{a})\geq\mu(x)
\]
for every $x\in G$.

Therefore
\[
\arraycolsep=.5mm\begin{array}{rl}
\mu(f(x_1^n))&=\mu(x_1\circ\varphi(x_2)\circ\varphi^2(x_3)\circ\ldots\circ\varphi^{n-1}(x_n)\circ
b)\\[4pt]
&\geq\min\{\mu(x_1),\mu(\varphi(x_2)),\mu(\varphi^2(x_3)),\ldots,\mu(\varphi^{n-1}(x_n)),\mu(b)\}
\\[4pt]
&\geq\min\{\mu(x_1),\mu(x_2),\mu(x_3),\ldots,\mu(x_n),\mu(b)\}\\[4pt]
&\geq\min\{\mu(x_1),\mu(x_2),\mu(x_3),\ldots,\mu(x_n)\},
\end{array}
\]
which proves that the first condition of the Definition~\ref{D2.3}
is satisfied.

To prove the second condition observe that from (\ref{skew}) and
(\ref{ho-2}) it follows
\[
\overline{x}=\left(\varphi(x)\circ\varphi^2(x)\circ\ldots\circ\varphi^{n-2}(x)\circ
b\right)^{-1}.
\]
Thus
\[
\arraycolsep=.5mm\begin{array}{rl}
\mu(\overline{x})&=\mu\left(\left(\varphi(x)\circ\varphi^2(x)\circ\ldots\circ\varphi^{n-2}(x)\circ
b\right)^{-1}\right)\\[4pt]
&\geq\mu\left(\varphi(x)\circ\varphi^2(x)\circ\ldots\circ\varphi^{n-2}(x)\circ
b\right)\\[4pt]
&\geq\min\{\mu(\varphi(x)),\mu(\varphi^2(x)),\ldots,\mu(\varphi^{n-2}(x)),\mu(b)\}
\\[4pt]
&\geq\min\{\mu(x),\mu(b)\}=\mu(x),
\end{array}
\]
which completes the proof.
\end{proof}
\begin{corollary}\label{C2.13}
If $(G,f)$ is a ternary group, then any fuzzy subgroup of
$ret_a(G,f)$ is a fuzzy ternary subgroup of $(G,f)$.
\end{corollary}
\begin{proof}
Since $\overline{a}$ is a neutral element of a group $ret_a(G,f)$
then $\mu(\overline{a})\geq\mu(x)$ for all $x\in G$. Thus
$\mu(\overline{a})\geq\mu(a)$. But in ternary group
$\overline{\overline{a}}=a$ for any $a\in G$, whence
$\mu(a)=\mu(\overline{\overline{a}})\geq\mu(\overline{a})\geq\mu(a)$.
So, $\mu(a)=\mu(\overline{a})\geq\mu(x)$ for all $x\in G$. This
means that the assumptions of Theorem~\ref{T2.12} are satisfied.
\end{proof}

\begin{example}\label{E2.14}
Consider the ternary group $(G,f)$, derived from the additive group
$\mathbb{Z}_4$. Let $\mu$ be a fuzzy subgroup of the group
$ret_1(G,f)$ induced by subgroups $S_1=\{11\}$, $S_2=\{5,11\}$ and
$S_3=\{1,3,5,7,9,11\}$, i.e., let $\mu(11)=t_1$, $\mu(5)=t_2$,
$\mu(1)=\mu(3)=\mu(7)=\mu(9)=t_3$ and $\mu(x)=t_4$ for $x\not\in
S_3$, where $0\leq t_4<t_3<t_2<t_1\leq 1$. Then
$\mu(\overline{5})=\mu(7)=t_3<t_2=\mu(5)$, which means that $\mu$ is
not a fuzzy ternary subgroup of $\mu(11)=t_1$.
\end{example}

From the above example it follows that:

\begin{enumerate}
\item There are fuzzy subgroups of $ret_a(G,f)$ which are not
fuzzy $n$-ary subgroup of $(G,f)$.
\item In Theorem~\ref{T2.12} the assumption $\mu(a)\geq\mu(x)$ cannot be omitted. In the
above example we have $\mu(1)=t_3<t_2=\mu(5)$.
\item The assumption $\mu(a)\geq\mu(x)$ cannot be replaced by the
natural assumption $\mu(\overline{a})\geq\mu(x)$ ($\overline{a}$
is the identity of $ret_a(G,f)$). In the above example
$\overline{1}=11$ and $\mu(11)\geq\mu(x)$ for all
$x\in\mathbb{Z}_{12}$ but $\mu$ is not a fuzzy $n$-ary subgroup of
$(\mathbb{Z}_{12},f)$.
\end{enumerate}

\begin{theorem}
Let $(G,f)$ be an $n$-ary group $b$-derived from the group
$(G,\circ)$. Any fuzzy subgroup $\mu$ of $(G,\circ)$ such that
$\mu(b)\geq\mu(x)$ for every $x\in G$ is a fuzzy $n$-ary group of
$(G,f)$.
\end{theorem}
\begin{proof}
The first condition of the Definition~\ref{D2.3} is obvious. To
prove the second observe that in an $n$-ary group $(G,f)$
$b$-derived from the group $(G,\circ)$
 \[
 \overline{x}=(x^{n-2}\circ b)^{-1},
 \]
where $x^{n-2}$ is the power of $x$ in $(G,\circ)$.

Therefore
 \[\arraycolsep=.5mm\begin{array}{rl}
\mu(\overline{x})&=\mu((x^{n-2}\circ b)^{-1})\geq \mu(x^{n-2}\circ
b)\geq\min\{\mu(x),\mu(b)\}=\mu(x)
\end{array}
 \]
for all $x\in G$. The proof is complete.
\end{proof}

\begin{corollary}\label{C2.14}
Any fuzzy subgroup of a group $(G,\circ)$ is a fuzzy $n$-ary
subgroup of an $n$-ary group derived from $(G,\circ)$.
\end{corollary}
\begin{proof} If an $n$-ary group $(G,f)$ is derived from the
group $(G,\circ)$ then $b=e$ and $\mu(e)\geq\mu(x)$ for all $x\in
G$.
 \end{proof}

\section{Characterizations of fuzzy $n$-ary subgroups}

\begin{lemma}\label{L3.1}
Two level subsets $\mu_{s},\mu_{t}$ $(s<t)$ of a fuzzy $n$-ary
subgroup $\mu$ of $G$ are equal if and only if there is no $x\in
G$ such that $s\leqslant\mu(x)<t$.
\end{lemma}
\begin{proof} Let $\,\mu_{s}=\mu_{t}\,$ for some $\,s<t$. If
there exists $\,x\in G\,$ such that $\,s\leqslant\mu(x)<t$, then
$\,\mu_{t}$ is a proper subset of $\,\mu_{s}$, which is a
contradiction. Conversely assume that there is no $\,x\in G\,$
such that $\,s\leqslant\mu(x)<t$. If $\,x\in\mu_{s}$, then
$\,\mu(x)\geqslant s$, and so $\,\mu(x)\geqslant t$, because
$\,\mu(x)\,$ does not lie between $\,s\,$ and $\,t$. Thus
$\,x\in\mu_{t}$, which gives $\mu_{s}\subseteq\mu_{t}$. The
converse inclusion is obvious since $\,s<t$. Therefore
$\mu_{s}=\mu_{t}$.
\end{proof}

\begin{proposition}\label{P3.2} Let $\mu$ and $\lambda$ be two fuzzy
$n$-ary subgroups of $G$ with the same family of levels. If\ {\rm
Im}$(\mu)=\{t_1,\ldots ,t_m\}$ and\ {\rm
Im}$(\lambda)=\{s_1,\ldots ,s_p\}$, where $t_1>t_2>\ldots >t_m$
and  $s_1>s_2>\ldots
>s_p$, then
\begin{itemize}
 \item[(i)] $m=p$,
 \item[(ii)] $\mu_{t_i}=\lambda_{s_i}$ for $i=1,\ldots , m$,
 \item[(iii)] if $\mu(x)=t_i$, then $\lambda(x)=s_i$ for $x\in G$ and $i=1,\ldots,m$.
\end{itemize}
\end{proposition}
\begin{proof} $(i)$ and $(ii)$ are obvious. To prove $(iii)$
consider $x\in G$ such that $\mu(x)=t_{i}$. If $\lambda(x)=s_{j}$
then $s_{j}\geq s_{i}$, i.e.,
$\lambda_{s_{j}}\subseteq\lambda_{s_{i}}$. Since
$x\in\lambda_{s_{j}}=\mu_{t_{j}}$, we obtain
$t_{i}=\mu(x)\geqslant t_{j}$, which gives
$\mu_{t_{i}}\subseteq\mu_{t_{j}}$. Consequently,
$\lambda_{s_{i}}=\mu_{t_{i}}\subseteq\mu_{t_{j}}=\lambda_{s_{j}}.$
Thus $\lambda_{s_{i}}=\lambda_{s_{j}}$. Lemma~\ref{L3.1} completes
the proof.
 \end{proof}

\begin{theorem}\label{T3.3}
Let $\mu$ and $\lambda$ be two fuzzy $n$-ary subgroups of $G$ with
the same family of levels. Then $\mu=\lambda$ if and only if \
{\rm Im}$(\mu)={\rm Im}(\lambda).$
\end{theorem}
\begin{proof} Let Im$(\mu)={\rm Im}(\lambda)=\{s_{1},...,s_{n}\}$
and $s_{1}>....s_{n}.$ By Proposition~\ref{P3.2} for each $x\in G$
there exists $s_{i}$ such that $\mu(x)=s_{i}=\lambda(x).$ Thus
$\mu(x)=\lambda(x)$ for all $x\in G$, which gives $\mu=\lambda$.
\end{proof}

\begin{theorem}\label{T3.4}
Let $\{H_i\,|\,i\in I\}$, where $I\subseteq [0,1]$, be a
collection of $n$-ary subgroups of $G$ such that
\begin{itemize}
 \item[(i)] $G=\bigcup_{i\in I}H_i$,
 \item[(ii)] $i>j\Longleftrightarrow H_i\subset H_j$ for all $i,j\in I$.
\end{itemize}
Then $\mu$ defined by $\mu(x)=\sup\{i\in I\,|\,x\in H_i\}$ is a
fuzzy $n$-ary subgroup of $G$.
\end{theorem}
\begin{proof} By Theorem~\ref{T2.7}, it is sufficient to show that every
nonempty level $\mu_k$ is an $n$-ary subgroup of $G$. Let $\mu_k$
be non-empty for some fixed $k\in [0,1]$. Then
\[
\begin{array}{l}
k=\sup\{i\in I\,|\,i<k\}=\sup\{i\in I\,|\,H_k\subset H_i\}\\
\text {or}\\
k\not =\sup\{i\in I\,|\,i<k\}=\sup\{i\in I\,|\,H_k\subset H_i\}.
\end{array}
\]
In the first case we have $\mu_k=\bigcap_{i<k}H_i$,because
\[
x\in\mu_k \ \Longleftrightarrow \ x\in H_i \ \text {for \ all} \
i<k \ \Longleftrightarrow \ x\in \bigcap_{i<k}H_k.
\]
In the second case, there exists $\varepsilon >0$ such that
$(k-\varepsilon,i)\cap I=\emptyset$. In this case
$\mu_k=\bigcup_{i\geq k}H_i$. Indeed. if $x\in\bigcup_{i\geq k}
H_i$, then $x\in H_i$ for some $i\geq k$, which gives $\mu(x)\geq
i\geq k$. Thus $x\in\mu_k$, i.e., $\bigcup_{i\geq k} H_i\subseteq
\mu_k$.

Conversely, if $x\not\in\bigcup_{i\geq k}H_i$, then $x\not\in H_i$
for all $i\geq k$, which implies $x \not\in H_i$ for all
$i>k-\varepsilon$, i.e., if $x\in H_i$ then $i\leq k-\varepsilon$.
Thus $\mu(x)\leq k-\varepsilon$. Therefore $x\not\in\mu_k$. Hence
$\mu_k\subseteq\bigcup_{i\geq k}H_i$, and in the consequence
$\mu_k=\bigcup_{i\geq k}H_i$. This completes the proof.
\end{proof}

\begin{theorem}\label{T3.5}  Let $\mu$ be a fuzzy set in $G$ and
let \ {\rm Im}$(\mu)=\{t_{0},t_{1},...,t_{m}\}$, where
$t_{0}>t_{1}>...>t_{m}$. If $H_{0}\subset H_{1}\subset\ldots
\subset H_{m}=G$ are $n$-ary subgroups of $G$ such that
$\mu(H_{k}\setminus H_{k-1})=t_{k}$ for $k=0,1,\ldots,m$, where
$H_{-1}=\emptyset$, then $\mu$ is a fuzzy $n$-ary subgroup.
\end{theorem}
\begin{proof}
For any fixed elements $x_1,\ldots,x_n\in G$ there exists only one
$k=0,1,...,m$ such that $f(x_1^n)$ belongs to $H_{k}\setminus
H_{k-1}$. If all $x_1,\ldots,x_n$ belongs to $H_k$, then at least
one lies in $H_{k}\setminus H_{k-1}$ because in the opposite case
$x_1^n\in H_{k-1}$ implies $f(x_1^n)\in H_{k-1}$ which is a
contradiction. So, in this case
 \[
\mu(f(x_1^n))=t_k=\min\{\mu(x_1),\dots,\mu(x_n)\} .
 \]
If $x_1,\ldots,x_n$ are not in $H_k$, then at least one of them
belongs to some $H_p\setminus H_{p-1}$, where $p>k$. Then
 \[
 \mu(f(x_1^n))=t_k\geq t_p\geq\min\{\mu(x_1),\dots,\mu(x_n)\}.
 \]
This means that the first condition of the Definition~\ref{D2.3}
is satisfied in any case.

The condition is obvious since $x\in H_k\setminus H_{k-1}$ implies
$\overline{x}\in H_k$. Thus $\mu(\overline{x})\geq\mu(x)$.
 \end{proof}

\begin{corollary}\label{C3.6} Let $\mu$ be a fuzzy set in $G$ with
{\rm Im}$(\mu)=\{t_{0},t_{1},\ldots,t_{m}\}$, where
$t_{0}>t_{1}>\ldots>t_{m}$. If \ $H_{0}\subset H_{1}\subset\ldots
\subset H_{m}=G\,$ are $n$-ary subgroups of $G$ such that
$\mu(H_{k})\geq t_{k}\,$ for $k=0,1,...,m$, then $\mu$ is a fuzzy
$n$-ary subgroup in $G$.
 \end{corollary}

\begin{corollary}\label{C3.7}
If \ {\rm Im}$(\mu)=\{t_{0},t_{1},\ldots,t_{m}\}$, where
$t_{0}>t_{1}>\ldots >t_{m}$, is the image of a fuzzy $n$-ary
subgroup $\mu$ in $G$, then all level subsets  $\mu_{t_{k}}$ are
$n$-ary subgroups of $G$ such that $\mu(\mu_{t_{0}})=t_{0}\,$ and
$\mu(\mu_{t_{k}}\!\!\setminus\mu_{t_{k-1}})=t_{k}\,$ for
$k=1,2,...,m$.
\end{corollary}
\begin{proof} By Theorem~\ref{T2.7} all level subsets $\mu_{t_{k}}$ are $n$-ary subgroups.
Clearly $\mu(\mu_{t_{0}})=t_{0}$. Since $\mu(\mu_{t_{1}})\geq
t_{1}$, then $\mu(x)=t_{0}\,$ for $\,x\in\mu_{t_{0}}\,$ and
$\mu(x)=t_{1}\,$ for $x\in\mu_{t_{0}}\!\!\setminus\mu_{t_{1}}$.
Repeating this procedure, we conclude that
$\mu(\mu_{t_{k}}\!\!\setminus\mu_{t_{k-1}})=t_{k}$ for
$\,k=1,2,\ldots,m$.
 \end{proof}

\begin{theorem}\label{T3.8}
Let $(G,f)$ be a unipotent $n$-ary group. If $\mu$ is a fuzzy
$n$-ary subgroup of $G$ with the image {\rm Im}$(\mu)=\{t_{i}:i\in
I\}$ and $\Omega =\{\mu_{t}:t\in {\rm Im}(\mu)\}$, then

\begin{enumerate}
\item[$(i)$] there exists a unique $t_{0}\in {\rm Im}(\mu)$ such that
$t_{0}\geq t$ for all $t\in {\rm Im}(\mu)$,
\item[$(ii)$] $G$ is the set-theoretic union of all
$\mu_{t}\in\Omega$,
\item[$(iii)$] the members of $\Omega$ form a chain,
\item[$(iv)$] $\Omega$ contains all level $n$-ary subgroups of $\mu$
if and only if $\mu$ attains its infimum on all $n$-ary subgroups
of $G$.
\end{enumerate}
\end{theorem}
\begin{proof} $(i)$ From the fact that in a unipotent
$n$-ary group $(G,f)$ there exists an element $\theta$ such that
$f(x_1^n)=\theta$ for all $x_1^n\in G$ it follows
 \[
t_0=\mu(c)=\mu(f(x_1^n))\geq\min\{\mu(x_1),\ldots,\mu(x_n)\} \]
for all $x_1^n\in G$. Whence we conclude $(i)$.

$(ii)$ If $x\in G$, then $t_x=\mu(x)\in {\rm Im}(\mu)$. This
implies $x\in\mu_{t_x}\subseteq\bigcup\mu_{t}\subseteq G$, where
$t\in {\rm Im}(\mu)$, which proves $(ii)$.

$(iii)$ Since $\mu_{t_{i}}\subseteq\mu_{t_{j}}\Longleftrightarrow
t_{i}\geq t_{j}$ for $i,j\in I$, then the family $\Omega$ is
totally ordered by inclusion.

$(iv)$ Suppose that $\Omega$ contains all level $n$-ary subgroups
of $\mu$. Let $S$ be an $n$-ary subgroup of $G$. If $\mu$ is
constant on $S$, then we are done. Assume that $\mu$ is not
constant on $S$. We consider two cases: $(1)$ $S=G$ and $(2)$
$S\subset G$. For $S=G$ let $\beta=\inf\,{\rm Im}(\mu)$. Then
$\beta\leq t\in {\rm Im}(\mu)$, i.e. $\mu_{\beta}\supseteq\mu_{t}$
for all $t\in {\rm Im}(\mu)$. But $\mu_{0}=G\in\Omega$ because
$\Omega$ contains all level $n$-ary subgroups of $\mu$. Hence
there exists $t'\in {\rm Im}(\mu)$ such that $\mu_{t'}=G$. It
follows that $\mu_{\beta}\supset\mu_{t'}=G$ so that
$\mu_{\beta}=\mu_{t'}=G$ because every level $n$-ary subgroup of
$\mu$ is an $n$-ary subgroup of $G$.

Now it sufficient to show that $\beta=t'$. If $\beta<t'$, then
there exists $t''\in {\rm Im}(\mu)$ such that $\beta\le t''<t'$.
This implies $\mu_{t''}\supset\mu_{t'}=G$, which is a
contradiction. Therefore $\beta=t'\in {\rm Im}(\mu)$.

In the case $S\subset G$ we consider the fuzzy set $\mu_{S}$
defined by
$$
\mu_{S}(x)=\left\{\begin{array}{ccl}
\alpha&for&x\in S,\\
0&for&x\in G\setminus S.
\end{array}
\right.
$$
From the proof of Proposition~\ref{P2.5} it follows that $\mu_{S}$
is a fuzzy $n$-ary subgroup of $G$.

Let
 \[
J=\{i\in I: \mu(y)=t_{i}\;\;for\;some\;y\in S\}
 \]
and $\Omega_{S}=\{(\mu_{S})_{t_{i}}: i\in J\}$. Noticing that
$\Omega_{S}$ contains all level $n$-ary subgroups of $\mu_{S}$,
then there exists $x_{0}\in S$ such that
$\mu(x_{0})=\inf\{\mu_{S}(x)\,|\, x\in S\}$, which implies that
$\mu(x_{0})=\{\mu(x)\,|\,x\in S\}$. This proves that $\mu$ attains
its infimum on all $n$-ary subgroups of $G$.

To prove the converse let $\mu_{\alpha}$ be a nonempty level
subset of a fuzzy $n$-ary subgroup $\mu$. If $\alpha=t$ for some
$t\in {\rm Im}(\mu)$, then clearly $\mu_{\alpha}\in\Omega$. If
$\alpha\ne t$ for all $t\in {\rm Im}(\mu)$, then there does not
exist $x\in G$ such that $\mu(x)=\alpha$. But $\mu_{\alpha}$ is
nonempty, so there exists $x_0\in G$ such $\mu(x_0)>\alpha$. Let
$H=\{x\in G:\mu(x)>\alpha\}$. Because
$\mu(\overline{x})\geq\mu(x)$ \ $x\in H$ implies $\overline{x}\in
H$. Moreover for $x_1^n\in H$ we have
 \[
\mu(f(x_1^n))\geq\min\{\mu(x_1),\ldots,\mu(x_n)\}>\alpha .
 \]
Hence $H$ is an $n$-ary subgroup of $G$. By hypothesis, there
exists $y\in H$ such that $\mu(y)=\inf\{\mu(x)\,|\,x\in H\}$. But
$\mu(y)\in {\rm Im}(\mu)$ implies $\mu(y)=t'$ for some $t'\in {\rm
Im}(\mu)$. Hence $\inf\{\mu(x)\,|\,x\in H\}=t'>\alpha$. Note that
there does not exist $z\in G\,$ such that $\alpha\le\mu(z)<t'$.
This gives $\mu_{\alpha}=\mu_{t'}$. Hence $\mu_{\alpha}\in\Omega$.
Thus $\Omega$ contains all level $n$-ary subgroups of $\mu$.
\end{proof}

\begin{proposition}\label{P3.9} Let $(G,f)$ be an $n$-ary
group such that every descending chain of its $n$-ary subgroups
terminates at finite step. If $\mu$ is a fuzzy $n$-ary subgroup in
$G$ such that a sequence of elements of \ ${\rm Im}(\mu)$ is
strictly increasing, then $\mu$ has a finite number of values.
\end{proposition}
\begin{proof} Assume that Im$(\mu)$ is not finite.
Let $0\leq t_{1}<t_{2}<\ldots\leq 1$ be a strictly increasing
sequence of elements of Im$(\mu)$. Every level subset $\mu_{t_i}$
is an $n$-ary subgroup of $G$. For $x\in\mu_{t_i}$ we have
$\mu(x)\geq t_{i}>t_{i-1}$, which implies $x\in\mu_{t_{i-1}}$.
Thus $\mu_{t_i}\subseteq\mu_{t_{i-1}}$. But for $t_{i-1}\in {\rm
Im}(\mu)$ there exists $x_{i-1}\in G$ such that
$\mu(x_{i-1})=t_{i-1}$. This gives $x_{i-1}\in\mu_{t_{i-1}}$ and
$x_{i-1}\not\in\mu_{t_i}$. Hence $\mu_{t_i}\subset\mu_{t_{i-1}}$,
and so we obtain a strictly descending chain
$\mu_{t_1}\supset\mu_{t_2}\supset\mu_{t_3}\supset\ldots$ of
$n$-ary subgroups, which is not terminating. This contradiction
completes the proof.
 \end{proof}

\begin{proposition}\label{P3.10} If every fuzzy $n$-ary subgroup
of $G$ has the finite image, then every descending chain of
$n$-ary subgroup of $G$ terminates at finite step.
\end{proposition}
\begin{proof}
Suppose there exists a strictly descending chain\
 \[
G=S_{0}\supset S_{1}\supset S_{2}\supset\ldots
 \]
of $n$-ary subgroups of $G$ which does not terminate at finite
step. We prove that $\mu$ defined by
 \[
\mu(x)=\left\{\begin{array}{ccl}
\frac{k}{k+1}&for&x\in S_{k}\setminus S_{k+1},\\[1mm]
1&for&x\in\bigcap S_{k},
\end{array}
\right.
 \]
where $k=0,1,2,\ldots$, is a fuzzy $n$-ary subgroup with an
infinite number of values.

If $f(x_1^n)\in\bigcap S_k$, then obviously
$\mu(f(x_1^n))=1\geq\min\{\mu(x_1),\dots,\mu(x_n)\}$.

If $f(x_1^n)\not\in\bigcap S_k$, then $f(x_1^n)\in S_p\setminus
S_{p+1}$ for some $p\geq 0$. Since $x_1^n\in\bigcap S_k$ implies
$f(x_1^n)\in \bigcap S_k$, at least one of $x_1,\ldots,x_n$ is not
in $\bigcap S_k$. Let $S_m$ be the smallest $S_k$ containing all
these elements. For $m>p$ we have $f(x_1^n)\in S_m\subseteq
S_{p+1}$, which contradicts to the assumption on $f(x_1^n)$. So,
$m\leq p$ and consequently
\[
\mu(f(x_1^n))=\frac{p}{p+1}\geq\frac{m}{m+1}=\min\{\mu(x_1),\ldots,\mu(x_n)\}.
 \]
This proves that $\mu$ satisfies the first condition of the
Definition~\ref{D2.3}. The second condition is obvious.

Hence $\mu$ is a fuzzy $n$-ary subgroup with an infinite number of
different values. Obtained contradiction completes our proof.
 \end{proof}

\begin{proposition}\label{P3.11} Every ascending chain of $n$-ary
subgroups of an $n$-ary group $G$ terminates at finite step if and
only if the set of values of any fuzzy $n$-ary subgroup of $G$ is
a well-ordered subset of $[0,1]$.
\end{proposition}
\begin{proof} If the set of values of a fuzzy $n$-ary subgroup
$\mu$ is not well-ordered, then there exists a strictly decreasing
sequence $\{t_{n}\}$ such that $t_{n}=\mu(x_{n})$ for some
$x_{n}\in G$. But in this case $n$-ary subgroups $B_{n}=\{x\in G
\,|\,\mu(x)\geq t_{n}\}$ form a strictly ascending chain, which is
a contradiction.

To prove the converse suppose that there exist a strictly
ascending chain $A_{1}\subset A_{2}\subset A_{3}\subset\ldots$ of
$n$-ary subgroups. Then $S=\bigcup\limits_{n\in N}A_{n}$ is an
$n$-ary subgroup of $G$ and $\mu$ defined by
$$
\mu(x)=\left\{\begin{array}{cl}
0&for\;\,x\not\in S\, ,\\[1mm]
\frac{1}{k}&where\;\,k=\min\{n\in N\,|\,x\in A_{n}\}
\end{array}
\right.
$$
is a fuzzy set on $G$.

We prove that $\mu$ is a fuzzy $n$-ary subgroup. The case when one
of $x_1,\ldots,x_n$ is not in $S$ is obvious. If all these
elements are in $S$ then also $f(x_1^n)\in S$. Let $k,m$ be
smallest numbers such that $x_1^n\in A_k$ and $f(x_1^n)\in A_m$.
Then $k\geq m$ and
 \[
\mu(f(x_1^n))=\frac{1}{m}\geq\frac{1}{k}=\min\{\mu(x_1),\ldots,\mu(x_n)\}.
 \]
So, $\mu$ satisfies the first condition of the
Definition~\ref{D2.3}. The second condition is obvious.

This means that $\mu$ is a fuzzy $n$-ary subgroup. Since the chain
$A_{1}\subset A_{2}\subset A_{3}\subset\ldots$ is not terminating,
$\mu$ has a strictly descending sequence of values. Obtained
contradiction proves that the set of values of any fuzzy $n$-ary
subgroup is well-ordered. The proof is complete.
\end{proof}

Let $\varphi$ be any mapping from an $n$-ary group $G_1$ to an
$n$-ary group $G_2$, and $\mu$ and $\lambda$ be fuzzy sets in
$G_1$ and $G_2$ respectively. Then the {\it image} $\varphi(\mu)$
and {\it pre-image} $\varphi^{-1}(\lambda )$ of $\mu$ and
$\lambda$ respectively, are the fuzzy sets defined as follows:
\[
\varphi(\mu )(y)=\left\{\begin{array}{cl}\displaystyle
 \sup_{x\in\varphi^{-1}(y)}\{\mu(x)\} & \text {if} \ \varphi^{-1}(y)\not =\emptyset \\
 0 & \text {if} \ \varphi^{-1}(y)=\emptyset ,
\end{array}
\right.
\]
\[
\varphi^{-1}(\lambda)(x)=\lambda(\varphi(x))
\]
for all $x\in G_1$ and $y\in G_2$. If $\varphi$ is a homomorphism
then $\varphi(\mu )$ is called the {\it homomorphic image of $\mu$
under} $\varphi$.

\begin{proposition}\label{P3.12} Let $\varphi$ be any mapping from an $n$-ary group
$G_1$ to an $n$-ary group $G_2$, and let $\mu$ be any fuzzy
$n$-ary subgroup of $G_1$. Then for $t\in (0,1]$ we have
\[
\varphi(\mu)_t=\bigcap_{t>\varepsilon>0}\varphi(\mu_{t-\varepsilon}).
\]
\end{proposition}
\begin{proof} Suppose that $t\in (0,1]$ and $y=\varphi(x)\in
G_2$. If $y\in\varphi(\mu )_t$ then
$\varphi(\mu)(\varphi(x))=\sup\limits_{x\in
\varphi^{-1}\varphi(x)} \{\mu(x)\}\geq t$. Therefore for every
real number $\varepsilon
>0$ there exists $x_0\in\varphi^{-1}(y)$ such that
$\mu(x_0)>t-\varepsilon$. So that for every $\varepsilon
>0$, $y=\varphi(x_0)\in\varphi(\mu_{t-\varepsilon })$, and hence $
y\in \bigcap\limits_{t>\varepsilon
>0}\varphi(\mu_{t-\varepsilon})$. Conversely, let $y\in\bigcap\limits_{t>\varepsilon
>0}\varphi(\mu_{t-\varepsilon})$, then for each $\varepsilon >0$
we have $y\in\varphi(\mu_{t-\varepsilon})$ and so there exists
$x_0\in \mu_{t-\varepsilon}$ such that $y=\varphi(x_0)$. Therefore
for each $\varepsilon >0$ there exists $x_0\in \varphi^{-1}(y)$
and $\mu(x_0)\geq t-\varepsilon$. Hence $\varphi(\mu
)(y)=\sup\limits_{x_i \in \varphi^{-1} (y)} \{\mu(x_i)\}\geq
\sup\limits_{t>\varepsilon
>0}\{t-\varepsilon\}=t$. So $y\in\varphi(\mu )_t$, and this
completes the proof.
\end{proof}

\begin{theorem}\label{T3.13} Let $\varphi$ be any homomorphism from an $n$-ary
group $G_1$ to an $n$-ary group $G_2$, and let $\mu$ be any fuzzy
$n$-ary subgroup of $G_1$. Then the homomorphic image $\varphi(\mu
)$ is a fuzzy $n$-ary subgroup of $G_2$.
 \end{theorem}
\begin{proof} By Theorem~\ref{T2.7}, $\varphi(\mu )$ is a fuzzy
$n$-ary subgroup if each nonempty level subset $\varphi(\mu )_t$
is an $n$-ary subgroup of $G_2$. If $t=0$ then $\varphi(\mu
)_t=G_2$ and if $t\in (0,1]$ then by Proposition~\ref{P3.12},
$\varphi(\mu)_t=\bigcap\limits_{t>\varepsilon
>0}\varphi(\mu_{t-\varepsilon})$. So $\varphi(\mu_{t-\varepsilon } )$
is nonempty for each $t>\varepsilon
>0$. Thus $\mu_{t-\varepsilon}$ is a nonempty level subset
of $\mu$ and by Theorem~\ref{T2.7} is an $n$-ary subgroup of
$G_1$. So, the homomorphic image $\varphi(\mu_{t-\varepsilon})$ is
an $n$-ary subgroup of $G_2$. Hence $\varphi(\mu )_t$ being an
intersection of a family of $n$-ary subgroups is also an $n$-ary
subgroup of $G_2$.
\end{proof}

\begin{theorem}\label{T3.14} Let $\varphi$ be a surjection from an $n$-ary
group $G_1$ to an $n$-ary group $G_2$, and let $\mu$ be a fuzzy
$n$-ary subgroup of $G_1$ which has the sup-property. If
$\{\mu_{t_i} \,|\,i\in I\}$ is the collection of all level $n$-ary
subgroups of $\mu$, then $\{\varphi(\mu_{t_i} )\,|\,i\in I\}$ is
the collection of all level $n$-ary subgroups of $\varphi(\mu)$.
 \end{theorem}
\begin{proof} Let $t\in [0,1]$, then
\[
u\in\varphi(\mu ) \ \Longrightarrow \varphi(\mu )(u)\geq t \
\Longrightarrow \ \sup\{\mu(x)\,|\,x\in\varphi^{-1} (u)\}\geq t.
\]
Since $\mu$ has sup-property, this implies that $\mu(x_0)\geq t$
for some $x_0\in \varphi^{-1}(u)$ Then $x_0\in\mu_t$ and hence
$\varphi(x_0)=u\in\varphi(\mu_t)$. Therefore, we have
$\varphi(\mu)_t \subseteq\varphi(\mu_t)$. Now, if
$u\in\varphi(\mu_t)$ then $u=\varphi(x)$ for some $x\not\in\mu_t$
and hence
\[
\varphi(\mu)(u)=\sup\{\mu(z)\,|\,z\in\varphi^{-1}(u)\}=
\sup\{\mu(z)\,|\,z\in\varphi(z)=\varphi(x)\}\geq\mu(x)\geq t.
\]
Therefore $u\in\varphi(\mu)_t$ and hence
$\varphi(\mu_t)\subseteq\varphi(\mu)_t$. Thus we have
$\varphi(\mu)_t=\varphi(\mu_t)$ for every $t\in [0,1]$. In
particular, $\varphi(\mu)_{t_i}=\varphi(\mu_{t_i})$ for all $i\in
I$. Hence all subsets $\varphi(\mu_{t_i})$ are level $n$-ary
subgroups of $\varphi(\mu)$. Also these are the only level $n$-ary
subgroups of $\varphi(\mu)$.
 \end{proof}

The following example shows that surjectiveness of $\varphi$ in
Theorem~\ref{T3.14} is essential.

\begin{example}\label{E2.15}
Let $(\mathbb{Z}_2, f)$ and $(\mathbb{Z}_4,g)$ be two ternary
groups derived from the additive groups $\mathbb{Z}_2$ and
$\mathbb{Z}_4$, respectively. Define $\varphi :\mathbb{Z}_2
\longrightarrow\mathbb{Z}_4$ by $\varphi(x)=x$ for
$x\in\mathbb{Z}_2$. Then $\varphi$ is not a surjective
homomorphism. Define $\mu :\mathbb{Z}_2\longrightarrow [0,1]$ by
$\mu(0)=0.3$ and $\mu(1)=0.1$. Then $\mu$ is a fuzzy ternary
subgroup of $(\mathbb{Z}_2,f)$ having sup-property. The level
ternary subgroups of $\mu$ are $\mu_{0.3}=\{0\}$ and
$\mu_{0.1}=G_1$. Now, $\varphi(\mu)$ is defined by
\[
\varphi(\mu)(0)=0.3, \;\; \varphi(\mu)(1)=\varphi(\mu)(3)=0, \;\;
\varphi(\mu)(2)=0.1.
\]
Hence the level ternary subgroups of $\varphi(\mu)$ are $\{0\}$,
$\{0,2\}$ and $\{0,1,2,3\}$. Therefore $\{\varphi(\mu_{0.3}),
\,\varphi(\mu_{0.1})\}$ does not contain all level ternary
subgroups of $\varphi(\mu)$.
 \end{example}

\section{Normal fuzzy $n$-ary subgroups}

\begin{definition}\label{D4.1}
Let $\mu$ be a fuzzy set of $G$. An element $\theta\in G$ is
called {\it $\mu$-maximal } if $\mu(\theta)\geq\mu(x)$ for all
$x\in G$. A fuzzy set $\mu$ with the property $\mu(\theta)=1$ is
called {\it normal}.
\end{definition}

Any fuzzy set $\mu$ with finite image has a $\mu$-maximal element.
In $n$-ary groups derived from binary group $(G,\circ)$ the
identity of $(G,\circ)$ is a $\mu$-maximal element for any fuzzy
subgroup of $(G,\circ)$. In unipotent $n$-ary groups the element
$\theta=f(x_1^n)$ is $\mu$-maximal for all fuzzy $n$-ary
subgroups. Thus a fuzzy $n$-ary subgroup $\mu$ of a unipotent
$n$-ary group is normal if and only if $\mu(\theta)=1$. Obviously
a characteristic function $\chi_A$ of any $n$-ary subgroup $A$ of
$G$ is normal.

\begin{proposition}\label{P4.2}
Let $\theta$ be a $\mu$-maximal element of a fuzzy $n$-ary
subgroup of an $n$-ary group $(G,f)$. Then a fuzzy set $\mu^+$
defined by $\mu^+(x)=\mu(x)+1-\mu(\theta)$ for all $x\in G$, is a
normal fuzzy $n$-ary subgroup of $G$ which contains $\mu$.
\end{proposition}
\begin{proof} Indeed,
\[
\arraycolsep=.5mm\begin{array}{rl}
\mu^+(f(x_1^n))&=\mu(f(x_1^n))+1-\mu(\theta)
\geq\min\{\mu(x_1),\mu(x_2),\ldots,\mu(x_n)\}+1-\mu(\theta)\\[4pt]
&=\min\{\mu(x_1)+1-\mu(\theta),\mu(x_2)+1-\mu(\theta),\ldots,\mu(x_n)+1-\mu(\theta)\}\\[4pt]
&=\min\{\mu^+(x_1),\mu^+(x_2),\ldots,\mu^+(x_n)\}\\[4pt]
\end{array}
\]
and
\[
\mu^+(\overline{x})=\mu(\overline{x})+1-\mu(\theta)\geq\mu(x)+1-\mu(\theta)=\mu^+(x).
\]
Clearly $\mu^+$ is normal and $\mu\subseteq\mu^+$.
 \end{proof}

It is clear that in a unipotent $n$-ary group a fuzzy set $\mu$ is
normal if and only if $\mu^+=\mu$.

\begin{corollary}\label{C4.3} Let $\mu$ and $\mu^+$ be as in the
above Proposition. If there is $x\in G$ such that $\mu^+(x)=0$,
then $\mu(x)=0$.
\end{corollary}

\begin{proposition}\label{P4.4} If a fuzzy $n$-ary subgroup $\mu$ of
an $n$-ary group has a $\mu$-maximal element, then
$(\mu^+)^+=\mu^+$. Moreover if $\mu$ is normal, then
$(\mu^+)^+=\mu$.
 \end{proposition}
\begin{proof} Straightforward.
 \end{proof}

\begin{proposition}\label{P4.5} Let $\mu$ be a fuzzy $n$-ary subgroup of
an $n$-ary group $G$. If there exists a fuzzy $n$-ary subgroup
$\nu$ of $G$ such that $\nu^+\subseteq\mu$, then $\mu$ is normal.
\end{proposition}
\begin{proof} Indeed, for $\nu^+\subseteq \mu$ we have
$1=\nu^+(\theta)\leq\mu(\theta)$. Hence $\mu(\theta)=1$.
 \end{proof}

Denote by $\mathcal{N}(G)$ the set of all normal fuzzy $n$-ary
subgroups of $G$. Note that ${N}(G)$ is a poset under the set
inclusion.

\begin{proposition}\label{P4.6} Let $\mu$ be a non-constant
fuzzy $n$-ary subgroup of an $n$-ary group $G$. If $\mu$ is a
maximal element of $(\mathcal{N}(G), \subseteq)$, then $\mu$ takes
only the values $0$ and $1$.
 \end{proposition}
\begin{proof} Observe that $\mu(\theta)=1$ since $\mu$
is normal. Let $x\in G$ be such that $\mu(x)\ne 1$. We claim that
$\mu(x)=0$. If not, then there exists $a\in G$ such that
$0<\mu(a)<1$. Let $\nu$ be a fuzzy set in $G$ defined by
$\nu(x)={\frac{1}{2}}(\mu(x)+\mu(a))$ for all $x\in G$. Then
clearly $\nu$ is well-defined, and
\[
\nu(\overline{x})={\frac{1}{2}}(\mu(\overline{x})+\mu(a))\geq
\frac{1}{2}(\mu(x)+\mu(a))=\nu(x)
 \]
for all $x\in G$. Moreover, for all $x_1^n\in G$ we get
 \[
 \arraycolsep=.5mm\begin{array}{rl}
\nu(f(x_1^n))&={\frac{1}{2}}(\mu(f(x_1^n)+\mu(a))
\geq\frac{1}{2}(\min\{\mu(x_1),\mu(x_2),\ldots,\mu(x_n)\}+\mu(a))\\[3mm]
&=\min\{\frac{1}{2}(\mu(x_1)+\mu(a)),\frac{1}{2}(\mu(x_2)+\mu(a))\ldots,
\frac{1}{2}(\mu(x_n)+\mu(a))\}\\[3mm]
&=\min\{\nu(x_1),\nu(x_2),\ldots,\nu(x_n)\}.
\end{array}
 \]
Hence $\nu$ is a fuzzy $n$-ary subgroup of $G$. It follows from
Proposition~\ref{P4.2} that $\nu^+\in {N}(G)$ where $\nu^+$ is
defined by $\nu^+(x)=\nu(x)+1-\nu(\theta)$ for all $x\in G$.
Clearly $\nu^+(x)\geq\mu(x)$ for all $x\in G$. Note that
 \[\arraycolsep=.5mm\begin{array}{rl}
\nu^+(a)&=\nu(a)+1-\nu(\theta)=
{\frac{1}{2}}(\mu(a)+\mu(a))+1-{\frac{1}{2}}(\mu(\theta)+\mu(a))\\[1mm]
&={\frac{1}{2}}(\mu(a)+1)>\mu(a)
\end{array}
 \]
and $\nu^+(a)<1=\nu^+(\theta)$. Hence $\nu^+$ is non-constant, and
$\mu$ is not a  maximal element of ${N}(G)$. This is a
contradiction.
 \end{proof}

We construct a new fuzzy $n$-ary subgroup from old. Let $t>0$ be a
real number. If $\alpha \in [0,1]$, $\alpha^{t}$ shall mean the
positive root in case $t<1$. We define $\mu^{t} : G\to [0,1]$ by
$\mu^{t}(x)=(\mu(x))^{t}$ for all $x\in G$.

\begin{proposition}\label{P4.7} If \ $\mu$ is a fuzzy $n$-ary
subgroup of an $n$-ary group $G$, then so is $\mu^{t}$. Moreover,
if $\theta$ is $\mu$-maximal, then $G_{\mu^{t}}=G_{\mu}$, where
$G_{\mu}=\{x\in G\, |\, \mu(x)=\mu(\theta)\}$.
\end{proposition}
\begin{proof}  For any $x,x_1^n\in G$, we have
$\mu^{t}(\overline{x})=(\mu(\overline{x}))^{t}\geq
(\mu(x))^{t}=\mu^{t}(x)$ and
 \[
 \arraycolsep=.5mm\begin{array}{rl}
\mu^{t}(f(x_1^n)&=(\mu(f(x_1^n))^{t}\geq (\min\{\mu(x_1),\ldots,\mu(x_n)\})^{t}\\[3mm]
&=\min\{(\mu(x_1))^{t},\ldots, (\mu(x_n))^{t}\}
=\min\{\mu^{t}(x_1),\ldots,\mu^{t}(x_n)\}.
\end{array}
 \]
Hence $\mu^{t}$ is a fuzzy $n$-ary subgroup. Moreover
\[
\arraycolsep=.5mm\begin{array}{rl}
 G_{\mu}&=\{x\in G\,|\,
\mu(x)=\mu(\theta)\}=\{x\in G\,|\,
(\mu(x))^{t}=(\mu(\theta))^{t}\}\\[3mm]
&=\{x\in G \,|\, \mu^t(x)=\mu^t(\theta)\}=G_{\mu^t}\, .
\end{array}
\]
This completes the proof.
 \end{proof}

\begin{corollary}\label{C4.8} If $\mu\in\mathcal{N}(G)$, then so is
$\mu^{t}$. \end{corollary}

\begin{definition}\label{D4.9} A fuzzy set $\mu$ defined on $G$ is
called {\it maximal} if it is non-constant and $\mu^+$ is a
maximal element of the poset $(\mathcal{N}(G), \subseteq)$.
\end{definition}

\begin{proposition}\label{P4.10} If $\mu$ is a maximal fuzzy
$n$-ary subgroup of an $n$-ary group $G$, then

\begin{enumerate}
\item[$(i)$] $\;\mu$ is normal,
\item[$(ii)$] $\;\mu$ takes only the values $0$ and $1$,
\item[$(iii)$] $\;G_{\mu}$ is a maximal $n$-ary subgroup of $G$.
\end{enumerate}
\end{proposition}
\begin{proof} Let $\mu$ be a maximal fuzzy $n$-ary subgroup of $G$. Then $\mu^+$
is a non-constant maximal element of the poset $(\mathcal{N}(G),
\subseteq )$. It follows from Proposition~\ref{P4.6} that $\mu^+$
takes only the values $0$ and $1$. Note that $\mu^+(x)=1$ if and
only if $\mu(x)=\mu(\theta)$, and $\mu^+(x)=0$ if and only if
$\mu(x)=\mu(\theta)-1$. By Corollary~\ref{C4.3}, we have
$\mu(x)=0$, that is, $\mu(\theta)=1$. Hence $\mu$ is normal, and
clearly $\mu^+=\mu$. This proves $(i)$ and $(ii)$.

$(iii)$ $\,G_{\mu}$ is a proper $n$-ary subgroup because $\mu$ is
non-constant. Let $S$ be an $n$-ary subgroup of $G$ containing
$G_{\mu}$. Noticing that, for any subsets $A$ and $B$ of $G$,
$A\subseteq B$ if and only if $\mu_A\subseteq\mu_B$, then we
obtain $\mu=\mu_{G_{\mu}}\subseteq\mu_S$. Since $\mu$ and $\mu_S$
are normal and $\mu =\mu^+$ is a maximal element of
$\mathcal{N}(G)$, we have that either $\mu=\mu_S$ or $\mu_S={\bf
1}$ where ${\bf 1}:G\to [0,1]$ is a fuzzy set defined by ${\bf
1}(x)=1$ for all $\,x\in G$. The later case implies that $S=G$. If
$\mu=\mu_S$, then $G_{\mu}=G_{\mu_S} =S$. This proves that
$G_{\mu}$ is a maximal $n$-ary subgroup of $G$.
 \end{proof}

\begin{definition}\label{D4.11} A normal fuzzy $n$-ary subgroup $\mu$ of $G$
is called {\it completely normal} if there exists $x\in G$ such
that $\mu(x)=0$. The set of all completely normal fuzzy $n$-ary
subgroups of $G$ is denoted by $\mathcal{C}(G)$.
\end{definition}

It is clear that $\mathcal{C}(G)\subseteq\mathcal{N}(G)$. The
restriction of the partial ordering $\subseteq$ of
$\mathcal{N}(G)$ gives a partial ordering of $\mathcal{C}(G)$.

\begin{proposition}\label{P4.13} Any non-constant maximal element of \
$(\mathcal{N}(G),\subseteq)$ is also a maximal element of \
$(\mathcal{C}(G),\subseteq)$.
 \end{proposition}
\begin{proof} Let $\mu$ be a non-constant maximal element of
$({N}(G), \subseteq)$. By Proposition~\ref{P4.6}, $\mu$ takes only
the values $0$ and $1$, and so $\mu(\theta)=1$ and $\mu(x)=0$ for
some $x\in G$. Hence $\mu\in\mathcal{C}(G)$. Assume that there
exists $\nu\in\mathcal{C}(G)$ such that $\mu\subseteq\nu$.
Obviously $\mu\subseteq\nu$ also in $\mathcal{N}(G)$. Since $\mu$
is maximal in $({N}(G),\subseteq )$ and $\nu$ is non-constant,
therefore $\mu=\nu$. Thus $\mu$ is maximal element of
$(\mathcal{C}(G), \subseteq)$.
 \end{proof}

\begin{proposition}\label{P4.13.}
Maximal fuzzy $n$-ary subgroup is completely normal.
\end{proposition}
\begin{proof} Let $\mu$ be a maximal fuzzy $n$-ary subgroup. By
Proposition~\ref{P4.10} $\mu$ is normal and $\mu=\mu^+$ takes only
the values $0$ and $1$. Since $\mu$ is non-constant, it follows
that $\mu(\theta)=1$ and $\mu(x)=0$ for some $x\in G$, which
completes the proof.
 \end{proof}

\begin{proposition}\label{P4.14}
Let $\theta$ be a $\mu$-maximal element of a fuzzy $n$-ary
subgroup of an $n$-ary group $G$. If \
$\varphi:[0,\,\mu(\theta)]\rightarrow [0,1]$ is an increasing
function, then a fuzzy set $\mu_{\varphi}$ defined on $G$ by
$\mu_{\varphi}(x)=\varphi(\mu(x))$ is a fuzzy $n$-ary subgroup.
Moreover, if $\varphi(t)\geq t$ for all $t\leq\mu(\theta),$ then
$\mu\subseteq\mu_{\varphi}$.
\end{proposition}
\begin{proof}
Since $f$ is increasing, then for all $x,x_1^n\in G$ we have
\[
\mu_{\varphi}(\overline{x})=\varphi(\mu(\overline{x})\geq
\varphi(\mu{x})=\mu_{\varphi(x)}
\]
and
 \[
 \arraycolsep=.5mm\begin{array}{rl}
\mu_{\varphi}(f(x_1^n))=\varphi(\mu(f(x_1^n)))&\ge\varphi(\min\{\mu(x_1),\ldots,\mu(x_n)\})
\\[3mm]
&=\min\{\varphi(\mu(x_1)),\ldots,\varphi(\mu(x_n))\}\\[3mm]
&=\min\{\mu_{\varphi}(x),\ldots,\mu_{\varphi}(y)\}\, .
\end{array}
\]
This proves that $\mu_{\varphi}$ is a fuzzy $n$-ary subgroup.

If $\varphi(t)\geq t$ for all $t\leq\mu(\theta)$, then
$\mu(x)\leq\varphi(\mu(x))=\mu_{\varphi}(x)$ for all $x\in G$,
which implies $\mu\subseteq\mu_f$.
\end{proof}


\begin{thebibliography}{13}
\bibitem{Bhat} P. Bhattacharya, N. P. Mukherjee, {\it Fuzzy relations and fuzzy groups},
Inform. Sci. {\bf 36} (1985), $267-282.$
\bibitem{Dor} W. D{\"o}rnte W, {\it Untersuchungen {\"u}ber einen
verallgemeinerten Gruppenbegriff}, Math. Z. {\bf 29} (1928),
$1-19.$
\bibitem{Rem} W. A. Dudek, {\it Remarks on $n$-groups}, Demonstratio
Math. {\bf 13} (1980), $165-181.$
\bibitem{auto} W. A. Dudek, {\it Autodistributive $n$-groups},
Commentationes Math. Annales Soc. Math. Polonae, Prace
Matematyczne {\bf 23} (1983), $1-11.$
\bibitem{D90} W. A. Dudek, {\it On $n$-ary group with only one skew
element}, Radovi Matemati{\v{c}}ki (Sarajevo), {\bf 6} (1990),
$171-175.$
\bibitem{Uni} W. A. Dudek: {\it Unipotent $n$-ary groups}, Demonstratio Math. {\bf 24}
(1991), $75-81.$
\bibitem{QRS7} W. A. Dudek, {\it Fuzzification of $n$-ary groupoids},
Quasigroups and Related Systems {\bf 7} (2000), $45-66.$
\bibitem{QRS8} W. A. Dudek, {\it On some old and new problems in $n$-ary groups}, Quasigroups
and Related Systems {\bf 8} (2001), $15-36.$
\bibitem{QRS13} W. A. Dudek, {\it Intuitionistic fuzzy approach to
$n$-ary systems}, Quasigroups and Related Systems {\bf 13} (2005),
$213-228.$
\bibitem{DGG} W. A. Dudek, K. G{\l}azek and B. Gleichgewicht, {\it A note on
the axioms of $n$-groups}, Colloquia Math. Soc. J. Bolyai {\bf 29}
("Universal Algebra", Esztergom (Hungary) 1977), $195-202.$
(North-Holland, Amsterdam 1982.)
\bibitem{DM1} W. A. Dudek and J. Michalski, {\it On a generalization of Hossz\'u theorem},
Demonstratio Math. {\bf 15} (1982), $783-805.$
\bibitem{DM} W. A. Dudek and J. Michalski, {\it On retracts of polyadic groups},
Demonstratio Math. {\bf 17} (1984), $281-301.$
\bibitem{Busse} J. W. Grzymala-Busse, {\it Automorphisms of polyadic automata},
J. Assoc. Comput. Mach. {\bf 16} (1969), $208-219.$
\bibitem{Hos} M. Hossz\'u: {\it On the explicit form of
$n$-groups}, Publ. Math. {\bf 10} (1963), $88-92.$
\bibitem{kasner} E. Kasner, {\it An extension of the group concept}
(reported by L. G. Weld), Bull. Amer. Math. Soc. {\bf 10} (1904),
$290-291.$
\bibitem{Ker} R. Kerner, {\it Ternary algebraic structures and their
applications in physics}, Univ. P. and M. Curie, Paris 2000.
\bibitem{Nik} D. Nikshych and L. Vainerman, {\it Finite quantum groupoids and
their applications}, Univ. California, Los Angeles 2000.
\bibitem{Poj} A. P. Pojidaev, {\it Enveloping algebras of Fillipov algebras},
Comm. Algebra {\bf 31} (2003), $883-900.$
\bibitem{Post} E. L. Post, {\it Polyadic groups}, Trans. Amer. Math. Soc. {\bf
48} (1940), $208-350.$
\bibitem{Ro} A. Rosenfeld, {\it Fuzzy groups}, J. Math. Anal. Appl. {\bf 35} (1971),
$512-517.$
\bibitem{Vai} L. Vainerman and R. Kerner, {\it On special classes of
$n$-algebras}, J. Math. Phys. {\bf 37} (1996), $2553-2565.$

\end{thebibliography}
\end{document}